\documentclass[11pt]{article}

\usepackage{amsmath}

\begin{document}


\title{The Langevin function and 
truncated exponential distributions}

\author{Grant Keady\\
Department of Mathematics,
Curtin University, Bentley, Australia
}

\date{\today}

\maketitle


\begin{abstract}
\noindent \emph{
Let $K$ be a random variable following
a truncated exponential distribution.
Such distributions are described by a single parameter
here denoted by $\gamma$.
The determination of $\gamma$ by Maximum Likelihood methods
leads to a transcendental equation.
We note that this can be solved in terms of
the inverse Langevin function.
We develop approximations to this guided by
work of Suehrcke and McCormick.
}

\medskip


\medskip


\end{abstract}


\section{Introduction, applications}

\noindent{Preamble, not for journal publication.} 
This paper was in 2007 submitted to {\it Solar Energy} as it was a slight improvement on the
theory part of a paper in that journal.
As it lacked observational data, that journal suggested it wasn't appropriate as it stood.
(WRF simulations have been considered, but not implemented.)
In 2012 another journal proposed a Special Issue on
Special Functions and Probablity and Statistics Applications, but this proposal didn't
eventuate. The arXiv version here is the 2012 version of the paper.

\bigskip

The frequency distribution of insolation values is important for 
the prediction of the performance of systems involving solar heating.
The fractional time distribution
for the clearness $K$,
for $K_{\rm min}\le K \le K_{\rm max}$  can be modelled by
\begin{equation}
F(K)
=\frac{\exp(\gamma K_{\rm min})-\exp(\gamma K)}
{\exp(\gamma K_{\rm min})-\exp(\gamma K_{\rm max})} ;
\label{eq:fK}
\end{equation}
see \cite{BCR} .
(Of course, $F(K)=0$ for $K<K_{\rm min}$ and $F(K)=1$ for $K>K_{\rm max}$.)
This is the doubly-truncated
exponential distribution.
We remark that $\gamma$ could be negative.
Equation~(\ref{eq:fK})
appears in equation (3.7.4) of \cite{R} p72  and equation (1) of \cite{SM}.
$K_{\rm max}$ and $K_{\rm min}$
are the maximum and minimum of the daily clearness index
over the period of interest (and a common approximation is
that $K_{\rm min}=0.05$).
Simulations of solar water and air heaters are discussed in~\cite{R}.
Oher situations -- and one which the author encounted is described in~\cite{D} --
might also be simulated in like manner.

\smallskip

The truncated exponential distribution
has many other applications including
distributions of earthquakes, of forest-fire sizes, 
raindrop sizes, reliability modelling, etc.
See \cite{AB,BD,HD} for example
(noting that, in some applications,
${K_{\rm min}}$, ${K_{\rm max}}$ and $\gamma$ are
all unknowns).
The Langevin function, defined in \S\ref{sec:Langevin}, is easily recognized in formulae in these papers.
For example, equation (2.7), $h_2(\nu)=r$, in \cite{BD}, is
$L(\nu/2)=2r-1$.

\smallskip

If one supposes $K_{\rm max}$ and $K_{\rm min}$ values to be known, it remains
to estimate $\gamma$. To do this, first calculate
\begin{equation}
\mu_K(\gamma)\!\!
= \!\! \int_{K_{\rm min}}^{K_{\rm max}}
\frac{dF}{dK} K\, dK 
=\!\!\!
\frac{K_{\rm min}\exp(\gamma K_{\rm min})- K_{\rm max}\exp(\gamma K_{\rm max})}
{\exp(\gamma K_{\rm min})-\exp(\gamma K_{\rm max})}
-\frac{1}{\gamma}.
\label{eq:gK}
\end{equation}
The quantity $\gamma$ can be estimated by solving
${\bar K}=  \mu_K(\gamma)$, 
where  ${\bar K}$ is the mean of the observations.
See \cite{R} p72 equation (3.7.5) and \cite{SM} equation (2).
The $\gamma$ so found is the method of moments estimate
and, in fact, the maximum likelihood estimate of $\gamma$.

\smallskip

To shorten the formulae which follow, define
$\alpha$ (an average) 
and
$\delta$ (half the difference) 
by
$$
\alpha= \frac{1}{2}( K_{\rm max}+K_{\rm min}) ,\qquad
\delta= \frac{1}{2}( K_{\rm max}-K_{\rm min}) .
$$
Using $K_{\rm min}=\alpha-\delta$,
and $K_{\rm max}=\alpha+\delta$ in
the right-hand side of equation~(\ref{eq:gK}) 
\begin{equation}
\mu_K(\gamma)
= -\frac{1}{\gamma}+\alpha +
{\delta}\,{\rm coth}({\gamma\delta})
\label{eq:muK}
\end{equation}
so equation~(\ref{eq:gK}e) becomes
\begin{equation}
\frac{({\bar K}-\alpha)}{\delta}
= -\frac{1}{\gamma\delta} + {\rm coth}({\gamma\delta}) .
\label{eq:3e}  
\end{equation}

\section{The Langevin function}
\label{sec:Langevin}

The purpose of this genuinely elementary
note is to call attention to the fact that
the function occuring on the right-hand side of equation~(\ref{eq:3e}) 
is the Langevin function, and various consequences of this.
For example,
$\gamma$ can be expressed exactly in terms of the 
inverse Langevin function.
Both the Langevin function and its inverse have
been studied widely in other contexts:
see, for example, \cite{C,P,SO}.
The Langevin function is defined, with $x\ne 0$, by
\begin{equation}
L(x)
:= {\rm coth}(x)-\frac{1}{x}
= \frac{d}{dx}\log\left(\frac{\sinh(x)}{x}\right)
\label{eq:Ldef}
\end{equation}
and, for $x=0$, by $L(0)=0$.
The function $L$ is a monotonic increasing function mapping the real line
into the interval from $-1$ to $1$, and it satisfies
\begin{equation}
\frac{d L(x)}{dx}
= 1- L(x)^2-\frac{2L(x)}{x} .
\label{eq:Lderiv}
\end{equation}
The {\it inverse Langevin function}
is the function $L^{-1}$ inverse to $L$.
Of course, the computation of $L^{-1}$ remains a numerical
issue, just as the solution for $\gamma$ of
equation~(\ref{eq:gK}e) is.
However, with the knowledge of the name of the function,
the solar applications can share algorithms and code
with the other applications, if it is so desired.
Knowing the name of the function involved assists in
the search for numerical software.
Numerical codes, in C, are available for
the inverse Langevin function.

\smallskip

With 
$$ 
x
=  \gamma \delta
= \frac{1}{2} \gamma (K_{\rm max}-K_{\rm min})
,\qquad
y
= \frac{{(\bar K}- \alpha)}{\delta}
= \frac{{2\bar K}- (K_{\rm max}+K_{\rm min})}
{K_{\rm max}-K_{\rm min}}
$$ 
equation~(\ref{eq:3e}) 
can be written
$y = L(x)$. 
Consequently
\begin{equation}
\gamma
= \frac{1}{\delta}L^{-1}\bigl( y\bigr) .
\label{eq:gammaLinvsol}
\end{equation}


\smallskip

Suehrcke et al. \cite{SM} gave an approximation to $\gamma$ solving
equation~(\ref{eq:3e}) 
which amounts to approximating
$x=L^{-1}(y)$ by
$$ x\approx A \tan(\frac{\pi y}{2}) .$$
Suehrcke used a polynomial
expression for $A$, finding its coefficients in the best way 
to be useful over the range of values
of $K_{\rm min}$, $K_{\rm max}$ and ${\bar K}$
occuring in the clearness application.
Restricting to 
the one tangent evaluation, this can be used in
other ways to approximate $L^{-1}$.
Let $\tau=\tan(\pi y/2)=\tan(\pi L(x)/2)$.
By fitting for small $\tau$ and for large $\tau$
using Maple we find that
$$ x\approx (\frac{6}{\pi})\tau\left(
\frac{1+b_n\tau^2}{1+b_d\tau^2}\right),
\qquad {\rm where\ } \ 
 b_n=\frac{\pi^2}{12} b_d, \quad
b_d =
\frac{20\pi^2-144}{\pi^2(60-5\pi^2)}
\approx 0.508 .
$$
This approximates $L^{-1}$ with a relative error of about 0.3\%
over the entire range and even smaller relative errors can be obtained
by using higher degree polynomials in the numerator and denominator.


\smallskip

The Langevin function -- and nothing worse than this --
arises in the formulae 
for {\it all} the higher moments 
and for {\it all} the higher cumulants
in the truncated exponential distribution.
The moment generating function $M$, and
the cumulant generation function $C$, are defined by
$$M(s) := 
\int_{K_{\rm min}}^{K_{\rm max}} \exp(sK) 
\frac{dF}{dK} \, dK,\qquad 
C(s):=\log(M(s)) .
$$
One easily finds
$$C(s)
= s\alpha + 
\ \log\left(\frac{\sinh(x+{s\delta})}{x+{s\delta}}\right) -
\ \log\left(\frac{\sinh(x)}{x}\right) .
$$
Differentiating this gives
\begin{equation}
C'(s) =\frac{dC}{ds}
= \alpha + {\delta} L(x + {s\delta}) .
\label{eq:Cderiv}
\end{equation}
The fact that $C'(0)=\mu_K(\gamma)$ recovers 
equation~(\ref{eq:muK}).
Equations~(\ref{eq:Lderiv}) and ~(\ref{eq:Cderiv})
show that each higher derivative of $C$ can be expressed as
a polynomial in $L=L(x+s\delta)$,
$\delta$, 
and $1/(x+{s\delta})$.
The most important cumulant after the mean $\mu_K(\gamma)$ is
the second, the variance $\sigma_K^2(\gamma)$.
\begin{equation}
\sigma_K^2(\gamma)
= C''(0)
=\delta^2 (1-L^2 -\frac{2L}{\gamma\delta}) ,
\label{eq:sigma2K}
\end{equation}
where $L=L(\gamma\delta)=(\mu_K(\gamma)-\alpha)/\delta$.
(The right hand side can be expressed, equally well,
in terms of $\coth(x)$ as in terms of $L(x)$.)

\smallskip

There is a practical use for equation~(\ref{eq:sigma2K}).
Let the first and second cumulants from the data
be denoted
${\bar K}$ (for the sample mean)
and $s_K^2$ (for the sample variance).
Next, denote by $v$ the value of 
the right-hand side of 
equation~(\ref{eq:sigma2K}),
when $\gamma$ is given by that solving 
${\bar K}=  \mu_K(\gamma)$. 
and $\mu_K(\gamma)$ is replaced by ${\bar K}$.
would give some indication of how
well the truncated exponential distribution (with the given
$K_{\rm min}$, $K_{\rm max}$) fits the data.

\smallskip

There is a more general
context for some of this, for example, in the context of natural exponential
families with quadratic variance functions.
Equation~(\ref{eq:sigma2K}) is an example of a
`quadratic variance function'.
The untruncated versions of these distributions are treated in~\cite{Mo},
but further consideration of such matters takes us too far from
the Langevin function observation which is the subject of this note.


\section{Conclusions}

We hope that practitioners in the different disciplines
using truncated exponential distributions will be able to benefit
from sharing experiences.
In particular, in
this paper, we have called attention to
the Langevin function, and the existing methods of 
approximating its inverse appropriate to estimating $\gamma$.


\end{document}